\documentclass[10pt]{article}
\usepackage{graphicx}
\usepackage{amsmath}
\usepackage{amssymb}
\usepackage{url}
\usepackage{bm}
\usepackage{xcolor}

\begin{document}
	
	\title{What is worthy of investigation? Philosophical attitudes and their impact on mathematical development by the example of discovering 10-adic numbers}
	
	\author{Merlin Carl\\{\small Europa-Universit\"at Flensburg}\\{\small Abteilung f\"ur Mathematik}\\{\small Auf dem Campus 1b}\\{\small 24943 Flensburg}\\{\small merlin.carl@uni-flensburg.de} \and Michael Schmitz\\{\small Europa-Universit\"at Flensburg}\\{\small Abteilung f\"ur Mathematik}\\{\small Auf dem Campus 1b}\\{\small 24943 Flensburg}\\{\small michael.schmitz@uni-flensburg.de}}
	\date{}
	
	\maketitle
	
	\begin{abstract}
	
	We describe in dialogue form a possible way of discovering and investigating 10-adic numbers starting from the naive question about a `largest natural number'. Among the topics we pursue are possibilities of extensions to transfinite 10-adic numbers, 10-adic representations of rational numbers, zero divisors, square roots and 10-adic roots of higher degree of natural numbers, and applications of 10-adic number representation in computer arithmetic. The participants of the dialogue are idealized embodiments of different philosophical attitudes towards mathematics. The article aims at illustrating how these attitudes interact, in both jarring and stimulating ways, and how they impact mathematical development.
		
%
		
	\end{abstract}
	
	\clearpage
	
\section{Preamble}
The following dialogue is an attempt to illustrate philosophical attitudes towards mathematics by the example of discovering 10-adic numbers. We want to illuminate how different philosophical attitudes interact, how they get in each other's way, but how they can also stimulate each other. Moreover, we want to exemplify how they impact the course of mathematical investigations on unknown territory. The idea for this came about from a conversation the second author had on a mathematics-enrichment camp in Kyoto with a senior high-school student, Yuki Kitamura, who was investigating 10-adic roots of natural numbers (within a school project) without prior knowledge about the theoretical background and without using common terminology. The discussion between the author and the student evolved from mutual misunderstanding in the beginning to an enriching experience for both sides in the end.\footnote{Using a computer Kitamura discovered the possible forms, $n \equiv 1 \pmod{40}$ or $n \equiv 9 \pmod{40}$, of natural numbers $n$ with $\gcd(10,n)=1$ that have 10-adic square roots, which will be derived in Section~\ref{sec:exploring} of this article.} The different philosophical attitudes, which actually came to light within only two people, are embodied here by five characters, each of whom idealizes a particular stance. More specifically, we have the following characters:

\begin{itemize}
	\item \textbf{Nana} represents a naive attitude. She does not question whether the expressions she writes are meaningful. She happily constructs and operates with concepts without worrying whether her objects are well-defined or the conditions for her operations are satisfied. For example, she has the opinion that every sequence of digits is a number. Although this might seem a bit careless, it turns out to be the initial spark a driving force of the whole investigation. 
	\item \textbf{Constance} is of a conservative persuasion. She is attached to the traditional scope of subject matter and relates expressions to it without considering alternative interpretations or deviations. Consequently, she argues that infinite numbers are no numbers, taking it for granted that such numbers must be natural.
	\item \textbf{Reila} adheres to critical realism. She insists that symbols must mean something, but is open to the possibility that they mean something other than originally intended or that one must first seek a meaning to a given system of symbols. Reila, however, would not consider a purely formal system sufficient in itself.
	\item \textbf{Forest}, being a formalist, lets the symbols stand for themselves and has no need for an object reference. Once the rules of manipulation are established, he is satisfied -- one must be able to deal with the expressions, then for him the meaning consists in the rules for their use.
	\item \textbf{Preston} approaches mathematics from a practitioner's perspective. He understands mathematical theories as tools that should be suitable to be beneficially applied in other field, inside or outside mathematics.
\end{itemize}

The dialogue is inspired both by Imre Lakatos' `Proofs and Refutations' \cite{Lakatos} and by Detlef Spalt's `Vom Mythos der Mathematischen Vernunft' (`On the Myth of Mathematical Reason', \cite{Spalt}), which also contains a mathematical investigation conducted by proponents of various directions in the philosophy of mathematics.

\section{The question and a first try}

Four mathematics students -- Constance, Forest, Reila, and Preston -- are sitting at the cafeteria together with a friend and his daughter Nana, a schoolchild.

\begin{description}
	
	\item[Friend:] You folks are smart at math, aren't you?
	
	\item[Reila:] Sure, why?
	
	\item[Friend:] The other day Nana asked me a mathematical question and I wasn't sure how to respond. I wonder what your answer would be.
	
	\item[Forest:] Shoot, Nana!
	
	\item[Nana:] What's the largest number?
	
	\item[Constance:] (smiles knowingly) See, Nana, there is no such number.
	
	\item[Nana:] Why's that? There must be one!
	
	\item[Constance:] Okay, then tell me the largest number you know!
	
	\item[Nana:] (ponders for a while) Thousand millions!
	
	\item[Constance:] Well, I would call that a billion, and it's not the largest number, because it's less then a billion and one.
	
	\item[Nana:] Then I say a thousand billions!
	
	\item[Constance:] (smiles complacently) That's a trillion, which is less than a trillion and one.
	
	\item[Nana:] Then I say a trillion trillions!
	
	\item[Constance:] Well, whatever that number's name is, I say a trillion trillions and one. Look Nana, no matter which number you tell me, it can never be the largest number, because I can always add one. And that's the reason why there is no such number.
	
	\item[Nana:] (looks grimly and folds her arms around her chest) You're playing unfair!
	
	\item[Constance:] (looks slightly puzzled) Why? It's just a mathematical fact. There's nothing unfair about it!
	
	\item[Nana:] Of course it is! I told you right at the beginning that I don't know the largest number. Therefore you can always add one to all numbers that \textit{I} know. But if I knew the largest number, you couldn't. And it's unfair that you don't tell me, ... or maybe you just don't know it yourself!\\
	
	\textit{Reila chokes on her coffee with laughter. Constance turns his palms up and looks somewhat perplexed.}\\
	
	\item[Forest:] Let me try to explain it another way. Nana, this largest number you imagine, what would be its last digit?
	
	\item[Nana:] $9$.
	
	\item[Forest:] Fine, and the second last?
	
	\item[Nana:] Again $9$.
	
	\item[Forest:] And the third last?
	
	\item[Nana:] Of course again $9$. And the same for the next digit and so on. All digits should be $9$.
	
	\item[Forest:] Agreed. Now tell me, how many digits has your largest number? I mean, where does your 9-line stop?
	
	\item[Nana:] They don't stop, they just go on and on. The number must go on forever, as otherwise it wouldn't be the largest number.
	
	\item[Forest:] But a number cannot have infinitely many digits. Numbers can have thousands and millions of digits, but always finitely many.
	
	\item[Nana:] That's not true! Of course there are such numbers, like $0.3333\ldots$ and so on. And they have infinitely many digits. So, I can also write my number, see $...999$?
	
	\item[Forest:] Look Nana, you cannot know this yet, but there's something really different about these numbers with infinitely many digits \textit{after} the decimal point. Believe me, there can't be infinitely many digits \textit{before} the decimal point. (turns to the others) I can't explain convergence and divergence of series to her, can I?
	
	\item[Nana:] You're out of arguments and just because I'm a child doesn't mean that I'm wrong.
	
	\item[Constance:] I have another idea how to explain what's wrong about your number. You know the addition algorithm?
	
	\item[Nana:] Of course I can add.
	
	\item[Constance:] Then would you please add one to your number. You can write it like this
	$$\begin{array}{rr}
		& ... 999999\\
		+	& 1\\\hline
		& \phantom{000000}\end{array}$$
	
	\textit{Nana complements the scheme in no time to}
	
	$$\begin{array}{rr}
		& ... 999999\\
		+	& 1\\\hline
		& ... 000000\end{array}$$
	
	\item[Forest:] And there you see, Nana.
	
	\item[Nana:] But that proves that I was right. The result is 0, so you don't get a larger number by adding 1. My number is the largest number.
	
	\item[Constance:] I don't agree. There is only one number that gives 0 when added to 1. And that number is $-1$. So your number would be $-1$, and that's not so large. Or your number is no number at all, and that would be my preferred conclusion.

	\item[Nana:] (Stares unhappily at the sheet of paper, shakes her head.) I think we did it wrong!
	
	\item[Reila:] What do you mean?
	
	\item[Nana:] The $...999+1=0$ thing. It's wrong. 
	
	\item[Constance:] (sighs) It's not, Nana. You applied the laws of arithmetic, and four mathematicians approved your calculation. 
	
	\item[Reila:] Nana, what do you think would be the ``right'' answer?
	
	\item[Nana:] We have $9+1=10$, and $99+1=100$, and $999+1=1000$ ...
	
	\item[Constance:] Oh no. I think I know what's coming next...
	
	\item[Nana:] ...and it goes on like that, the $9$s become $0$s and there is a new $1$...
	
	\item[Constance:] Please don't say what I think you want to say...
	
	\item[Nana:] So $...999+1$ must be $1....000$ -- a $1$ followed by infinitely many $0$s.
	
	\item[Constance:] Look, Nana. That's not sensible. You cannot have a $1$ followed by infinitely many $0$s.
	
	\item[Nana:] Why not? It makes a lot of sense. It was really strange that $...999$ should be such a small number, but now I understand that it is actually as huge as it should be.
	
	\item[Constance:] Because there is no position infinity in a number representation. There is a first, a second, a third, a hundreds, a billions. But these are all finite. None of them has infinitely much space to its right.  But no place to write your $1$!
	
	\item[Nana:] Well, now there is one. I just made it.
	
	\item[Constance:] (waves hands in desperation) Nana, that's absurd! I mean, where would you even write your $1$? There would be no place in the universe left for it once all of your $0$s are there!
	
	\item[Nana:] I could make the zeroes very small. 
	
	\item[Reila:] Even very small zeroes would not help much when you want them to be all of the same size. But if, for example, the $i$th zero was of size $2^{-i}$...
	
	\item[Forest:] Before Constance suffers a heart attack, let me try to get this a bit clearer: I will write $\overleftarrow{s}$ for infinitely many repetitions of a digit string $s$ to the left. Now Nana proposes $\overleftarrow{9}+1=1\overleftarrow{0}$. I wonder if that works out. Nana, what is $\overleftarrow{9}+\overleftarrow{1}$?
	
	\item[Nana:] That's easy:
	$$\begin{array}{rrr}
		& 0 & ...999999\\
		+ &	0 & ...111111\\\hline
		& 1 & ...111110\end{array}$$

	Which is what we should expect; after all, it is just $9\cdot\overleftarrow{1}+\overleftarrow{1}$, i.e., $10\cdot\overleftarrow{1}$, and multiplication by $10$ shifts everything by one place to the left and appends a $0$ on the right. 
	
	\item[Constance:] Would you please tell me which digit is ``shifted by one place to the left'' to yield your new leading digit with its infinitely many successors?
	
	\item[Nana:] Oh, come on, you know what I mean. You are just being stubborn. 
	
	\item[Constance:] No, in fact, I don't. 
	
	\item[Forest:] I believe Constance. The point is that she doesn't even want to know, due to her tendency to believe that everything that is worth knowing is already known. But she still has a point: A simple ``shifting'' does not account for the effect that Nana wants to achieve. 
	
	\item[Constance:] Nor will anything else. Forest, I respect your good will to extract meaning from even the most absurd proposals, but this is really crossing the border. There is just no reason to add this $1$, nor a place for it.
	
	\item[Nana:] But it has to be there! There are all of these carries, and we cannot simply ignore them. And if we can say that there are infinitely many digits, why can't we say that there is an infinite position?
	
	\item[Forest:] Constance, I respect your critical attitude insofar it triggers progress by forcing us to be precise, but I reject it when it is used to discourage research just because things seem unusual. Maybe Nana has something in her mind that makes sense; then it is worth investigating. It is of course always possible that things turn out to be inconsistent - in fact, your beloved traditional theories fare no better in this respect -, but that is the result of an investigation and not something we can see before we even try. 
	
	\item[Reila:] Right. We have to get the picture of what Nana is either seeing or believes to be seeing; then we need to find out whether it is actually a picture of something or a mere illusion. 
	
	\item[Forest:] Not at all. I distrust your so-called ``seeing''; it turns out to yield mere illusions too often to be more than a source of inspiration at best. The right way to proceed is to describe as precise as possible how Nana deals with her new objects, and then to consider whether the manipulation rules thus extracted 
	are consistent. 
	
	\item[Reila:] I disagree with you on the status of seeing things. Sure, there are misperceptions, but that does not prevent perception from being the way that things occur to us. On the other hand, your focus on formalisms may simply produce pedantic descriptions of nothing. But your method has helped as a vehicle of seeing things before, and I am willing to learn. So please go on. 
	
	\item[Forest:] Right; so back to business. Nana, that is an interesting point; so you propose to focus on the carry. Let me see. What is $\overleftarrow{89}+\overleftarrow{01}$?

	\item[Nana:] Do you really need to ask? Clearly, that is $\overleftarrow{90}$. 
	
	\item[Forest:] No leading digit here? After all, we get a carry $1$ infinitely often. 
	
	\item[Nana:] No, because they always go away afterwards. 
	
	\item[Forest:] I see. OK, in the addition of two such numbers, no carries greater than $1$ can occur. So you are saying: When the carry is eventually constantly $1$ when carrying out the addition from right to left, there is a leading $1$, and if it is $0$ unboundedly often, there isn't?
	
	\item[Nana:] (frowns) I guess.
	
	\item[Forest:] I see. 
	
	\item[Reila:] No, Forest - you do not ``see'' anything, if we believe in your story. That is a typical example of the weakness of your method, Forest. Now you have an addition rule, but still nothing to add with it. However, I think I can help now: What Nana seems to be considering are simply digit strings whose length 
	are transfinite ordinals. Are you folks familiar with those?\\
	
	\textit{Heads are shaken from right to left and back, hands are waved in typical ``more or less''-gestures.}\\
	
	\item[Reila:] Well, just read them up, then. It's fascinating stuff.\footnote{A good place to do so would be \cite{Devlin}, chapter $2$.} Roughly, an ordinal is just the set of its predecessors. So the first ordinal, which has no predecessors, is $\emptyset$, which we call $0$, then comes $\{\emptyset\}$, also known as ``$1$'', followed by $\{\emptyset,\{\emptyset\}\}$, which is $2$, and so on. Whenever we have such an initial segment of the ordinals, we regard it as a new set to get the next ordinal; thus, in particular, we obtain the first transfinite ordinal $\omega=\{0,1,2,3,...\}$; the numbers we talked about before are digit strings of length $\leq\omega$. But the ordinals go on after that: After $\omega$ comes $\omega\cup\{\omega\}=\{0,1,2,...,\omega\}$, known as $\omega+1$, $\omega+1\cup\{\omega+1\}$, which is $\omega+2$ etc. Now Nana simply proposes to consider digit strings of length $\omega+1$. 
	
	\item[Forest:] I see. That is indeed nice. But is it enough? What about, say $10\cdot 1\overleftarrow{0}$ or $10000\cdot 1\overleftarrow{0}$?
	
	\item[Reila:] Oh right. That should probably be $10\overleftarrow{0}$ and $10000\overleftarrow{0}$. So we should look at digit strings of length $<\omega\cdot 2$ to ensure closure under addition. 
	So we have a structure $(T,+)$, where $T$ is the set of functions from ordinals to $\{0,1,...,9\}$ and $+$ is defined like the before discussed addition on the $\omega$-part and like ordinary finite addition on the transfinite part, with the extra point that, if the carry in the finite part eventually becomes constantly $1$, 
	we add another $1$ to the transfinite part. 
	
	\item[Forest:] Right, that seems to be what Nana wants.\\
	
	\textit{Nana looks totally perplexed.}\\
	
	\item[Constance:] But mathematics is still not Santa Claus who gives you what you want. All you have is a set and a function. You call that weird function an ``addition'', but does it have any properties of addition? Is it commutative, associative? What about subtraction? Is there, for all $a$, $b$, a $c$ such that $a+c=b$ or $b+c=a$? Probably, you want to define subtraction in a similar way to addition by ``generalizing'' the usual algorithm to your transfinite monstrosities. But will this have the required property that $a+(b-a)=b$?
	
	\item[Forest:] Neither is mathematics a folklore festival. Terms and symbols have the meaning we give them by definition. The new addition is what we defined it to be, nothing else. You, Constance, may be used to the word ``addition'' as refering to an operation having certain properties, but your conditioning is not my concern. 
	
	\item[Reila:] I disagree. Given that this operation arose as and is intended to be an extension of addition, it is very natural to ask what of the properties of addition hold true in the new sense.  So these are good questions. But they are not very hard to answer.
	
	\item[Constance:] OK, I'll grant you commutativity, which is indeed obvious. After all, both $10$-adic and finite addition are commutative, and the carries also do not depend on the order of the summands. 
	
	\item[Reila:] Right; so let us consider associativity. Say we have three of these transfinite numbers, namely $a=a_{0}a_{1}$, $b=b_{0}b_{1}$, $c=c_{0}c_{1}$, where $a_{0}$, a (possibly empty) digit string is the transfinite initial part of $a$ and $a_{1}$ is the digit string of the last $\omega$ (or less) digits of $a$; and likewise for $b$ and $c$. Now, like in Constance's argument, associativity works for $a_{0}$, $b_{0}$ and $c_{0}$ and likewise for $a_{1}$, $b_{1}$ and $c_{1}$; so the last $\omega$ many digits of $a+(b+c)$ and $(a+b)+c$ will both agree with $(a_{1}+b_{1})+c_{1}$. 
	Moreover, we have $(a_{0}+b_{0})+c_{0}=a_{0}+(b_{0}+c_{0})$; so all that remains to be considered is the contribution of the carry. 
	
	\item[Constance:] Excellent idea. I propose to consider the carry for the case $a=\overleftarrow{69}$, $b=\overleftarrow{92}$ and $c=\overleftarrow{26}$. Then 
	$(a+b)$ is calculated by\\
	
	\begin{tabular}{ccccccccccc} 
		& & ... & $6_{\phantom{1}}$ & $9_{\phantom{1}}$ & $6_{\phantom{1}}$ & $9_{\phantom{1}}$ & $ 6_{\phantom{1}} $ & $ 9_{\phantom{1}} $ & $ 6_{\phantom{1}}$ & $ 9_{\phantom{1}} $ \\
		+ &  &... & $9_{1}$ & $2_{1}$ & $9_{1}$ & $2_{1}$ & $9_{1}$ & $2_{1}$ & $9_{1}$ & $2_{\phantom{1}}$ \\
		\hline
		& 1 & ... & $ 6_{\phantom{1}} $ &$ 2_{\phantom{1}} $ &$ 6_{\phantom{1}} $ &$ 2_{\phantom{1}} $ &$ 6_{\phantom{1}} $ &$ 2_{\phantom{1}} $ &$ 6_{\phantom{1}} $ &$ 1_{\phantom{1}} $ 
	\end{tabular}\\
	
	So $(a+b)+c$ gives us\\
	
	\begin{tabular}{ccccccccccc}
		& 1 & ... & 6 & 2 & 6 & 2 & 6 & 2 & 6 & 2 \\
		+ & & ...& 2 & 6 & 2 & 6 & 2 & 6 & 2 & 6\\
		\hline
		& 1 & ... & 8 & 8 & 8 & 8 & 8 & 8 & 8& 7
	\end{tabular}\\
	
	And this is $1\overleftarrow{8}7$. Nana, why don't you compute $a+(b+c)$?
	
	\item[Nana:] Why would I? It is clear that the result will be $1\overleftarrow{8}{7}$, as you just demonstrated. 
	
	\item[Reila:] No, Nana - we cannot just take associativity for granted, we are checking it. Let me do it: First, we obtain $b+c$:\\
	
	\begin{tabular}{cccccccc}
		& ... & $9_{\phantom{1}}$ & $2_{\phantom{1}}$ & $9_{\phantom{1}}$ & $2_{\phantom{1}}$ & $9_{\phantom{1}}$ & $2_{\phantom{1}}$\\
		+ & ... & $2_{\phantom{1}}$ & $6_1$ & $2_{\phantom{1}}$ & $6_1$ & $2_{\phantom{1}}$ & $6_{\phantom{1}}$\\
		\hline
		& ...& $1_{\phantom{1}}$ & $9_{\phantom{1}}$ & $1_{\phantom{1}}$ &  $9_{\phantom{1}}$ & $1_{\phantom{1}}$ & $8_{\phantom{1}}$              
	\end{tabular}\\
	
	So $b+c=\overleftarrow{91}8$.\\
	
	And next, $a+(b+c)$ yields:\\
	
	\begin{tabular}{cccccccc}
		&...& $6_{\phantom{1}}$ & $9_{\phantom{1}}$ & $6_{\phantom{1}}$ & $9_{\phantom{1}}$ & $6_{\phantom{1}}$ & $9_{\phantom{1}}$\\
		+&...& $1_{1}$ & $9_{\phantom{1}}$ & $1_{1}$ & $9_{\phantom{1}}$ & $1_{1}$ & $8_{\phantom{1}}$\\
		\hline
		&...& $8_{\phantom{1}}$ & $8_{\phantom{1}}$ & $8_{\phantom{1}}$ & $8_{\phantom{1}}$ & $8_{\phantom{1}}$ & $7_{\phantom{1}}$
	\end{tabular}\\
	
	So we get $\overleftarrow{8}7$. 
	
	\item[Constance:] ...which is clearly not the same as $1\overleftarrow{8}7$. So even associativity fails. Your ``new area'' crumbles to dust, as I predicted it would.
	
	\item[Nana:] What? That is quite confusing. I need to think about what went wrong.
	
	\item[Constance:] It is quite obvious what went wrong: We wrote down non-sensible sign combinations, and then found that they were, indeed, nonsensical. Hardly surprising. 
	
	\item[Forest:] We just got a non-associative addition. I still see nothing wrong with it.
	
	\item[Reila:] I do. But still, whether something makes sense is the result of an investigation, not something that should be judged by gut feeling in advance to prevent the investigation from happening. And in this case it just turned out ...
	
	\item[Nana:] OK, now I see what is wrong. The carry doesn't get away just because it gets absorbed from time to time.
	
	\item[The others:] What?
	
	\item[Nana:] Well, look: If we add $\overleftarrow{9}+1$, there is always a $1$, but it is not the same $1$: The first $1$ vanishes in the second position, the second in the third and so on. So, they all vanish, but still, the leftmost one remains...
	
	\item[Constance:] ...stubbornly stuck in your head.
	
	\item[Nana:] So, I think we have, for example, $\overleftarrow{08}+\overleftarrow{02}=1\overleftarrow{10}$, not $\overleftarrow{10}$, as I first thought. And then, for Constance's example, we get $1\overleftarrow{8}7$ for $a+(b+c)$, so all is well again.
	
	\item[Forest:] I think what Nana is proposing here is a new definition of addition, in which the roles of $0$ and $1$ are switched with respect to the carry, compared with our first attempt: If the carry is $1$ unboundedly often, we add $1$ to the transfinite part, and only when it is constantly $0$ we have a ``transfinite carry'' of $0$.

	\item[Constance:] And again, this is a non-starter. Just let $a=\overleftarrow{65}$, $b=\overleftarrow{56}$, $c=\overleftarrow{05}$. Then, by your new rules, we have $(a+b)+c=1\overleftarrow{72}6$ and $a+(b+c)=2\overleftarrow{72}6$.

\end{description}

\section{Constituting a new area}

\textit{After a short moment of silence...}

\begin{description}
	
	\item[Nana:] Then let's do it the way Constance suggested.
	
	\item[Constance:] What did I suggest?
	
	\item[Nana:] That $\overleftarrow{9}+1=0$.
	
	\item[Constance:] (eyes bulging) I didn't suggest that! I said that to convince you that your infinite 9-line is nonsensical.
	
	\item[Nana:] If I'm not mistaken, we then don't have the $(a+b)+c\ne a+(b+c)$ trouble any more. And that's what you wanted.
	
	\item[Constance:] Yes, Nana. The associate law should be valid if we calculate this way. But you forgot that we then have $\overleftarrow{9} = -1$.
	
	\item[Nana:] That doesn't bother me.
	
	\item[Constance:] Then tell me. What number is larger, $9$ or $-1$?
	
	\item[Nana:] $9$ is larger.
	
	\item[Constance:] And $99$ or $9$?
	
	\item[Nana:] $99$.
	
	\item[Constance:] And $999$ or $99$?
	
	\item[Nana:] I get it, Constance, the more $9$s we use, the larger gets the number.
	
	\item[Constance:] Then you should also understand that $\overleftarrow{9} = -1 < 1$ is contradictory.
	
	\item[Nana:] Of course it's weird. But didn't you say before, it's a big difference whether a number goes on forever or not? So, the new number $\overleftarrow{9}$ is somewhat strange, but still it exists.
	
	\item[Forest:] I think Nana has a point. Of course we can consider infinite sequences of digits, write them from right to left, and use the addition Nana proposes. And we can scrutinize what kind of mathematical structure we get.
	
	\item[Constance:] Correct me if I'm wrong. In mathematics, when a theory is extended, we have to do it in a way that ensures that the old rules still apply. And that's obviously not given in this case. So I don't see the point.
	
	\item[Reila:] I agree with Forest, and indeed I must correct you Constance. There are loads of examples where a theory is extended and the old rules have to be revised. For instance, in the naturals a product is never smaller than any of its factors. This rule neither holds for integers nor rationals. Does this imply that we shouldn't extend the additive semigroup $(\mathbb{N},+)$ to $(\mathbb{Z},+)$ and so on?
	
	\item[Constance:] (hesitates)
	
	\item[Forest:] Thanks for your support, Reila. Let's go ahead.
	
	\item[Reila:] So, we have a new representation for $-1$, namely $\overleftarrow{9}$. And when you think about it, it's quite consistent.
	
	\item[Constance:] Consistent with what?
	
	\item[Reila:] Nana, what do you think: Is $0.\overline{9} = 1$ or not?\footnote{The following line of argumentation was inspired by James Tanton's excellent video lecture series \cite{Tanton}.}
	
	\item[Nana:] I think it's almost true. Because $0.999...$ gets nearer and nearer to $1$, so it should be pretty close to 1.
	
	\item[Reila:] Let me convince you that it is actually equal to 1.
	
	\item[Constance:] Didn't you want to explain something else, Reila?
	
	\item[Reila:] Hang on, Constance. I will come to that point.
	
	\item[Reila:] Let's give this number a name, say $x=0.\overline{9}$. I claim that $x=1$, agree?
	
	\item[Nana:] Agreed.
	
	\item[Reila:] What's $10x$?
	
	\item[Nana:] It's $9.\overline{9}$.
	
	\item{Reila:} And what's $9x$ then?
	
	\item[Nana:] As $9x=10x-x$ that should be ... ah, I see! $9x=9$, so $x=1$. That's convincing!
	
	\item[Reila:] Is that convincing, Constance?
	
	\item[Constance:] Sure, I've seen that before.
	
	\item[Reila:] Okay, then back to the number we are interested in. It also needs a name, say $y=\overleftarrow{9}$. Nana, what's $10y$?
	
	\item[Nana:] That's $\overleftarrow{9}0$.
	
	\item[Reila:] Correct. Then tell me, what $-9y$ is? I give you a hint, do $y-10y$.
	
	\item[Nana:] Okay, so I have to compute $...999-...990$, and that's clearly $9$. And -- wow -- that means $-9y = 9$, so $y=-1$.
	
	\item[Reila:] Well done, Nana! See, Constance. That's what I meant by consistent.
	
	\item[Constance:] But you forgot to mention, that the rules that you apply are valid for representations with infinitely many digits after the decimal point as they are converging series. But they don't apply for these new bogus numbers. So it's still pointless.
	
	\item[Forest:] I strongly disagree. We consider sequences of digits, have agreed on a definition for addition, and I think all of us can imagine how multiplication should be defined. And I think that we all have enough overview to see that from these definitions the rules that are needed for Reila's trick should follow easily. And in my opinion, that's what mathematics is all about. To agree on certain prerequisites and scrutinize what they imply.\footnote{This is a typical implicationist view on mathematics as for instance described by Potter in \cite{Potter}. p.~9: \textit{`[...], the axiomatic method was by the 1920s becoming such a
			mathematical commonplace, and implicationism such a common attitude towards it, that it was inevitable it would be applied to the recently founded theory of sets. [...]
			One of the evident
			attractions of the implicationist view of set theory is that it obviates the tedious
			requirement imposed on the realist to justify the axioms as true and replaces it
			with at most the (presumably weaker) requirement to persuade the reader to
			be interested in their logical consequences. Even in the extreme case where
			our axiom system turned out to be inconsistent, this would at worst make its
			consequences uninteresting, but we could then convict the implicationist only
			of wasting our time, not of committing a mistake.'}
	}
	
	\item[Constance:] But one agrees only on prerequisites that are true, otherwise it's merely doing smart jokes.
	
	\item[Forest:] Remember the story where Russell deduced from the assumption that $0=1$ that he is the pope? I like that one. But seriously, I don't agree. Do parallel lines meet in a point at infinity? In Euclidean Geometry it's false, but in Projective Geometry it's true. And you wouldn't call one of these senseless, would you? The point is that we cannot judge beforehand whether a prerequisite is true or not. It's in the nature of prerequisites, or better say axioms, that we can't. 
	
	\item[Reila:] Anyway, I still agree with Constance that there should be a reasonable meaning for a theory to be worth studying. But there still could be a sensible interpretation for our new `infinite numbers'; maybe we just haven't found it yet. And of course, that meaning must be different from the usual meaning in the decimal system...
	
	\item[Preston:] And by meaning you mean some sort of application?
	
	\item[Reila:] Not necessarily...
	
	\item[Forest:] We will never know if we keep discussing and don't start analyzing.
	
	\item[Constance:] Then go ahead, for god's sake!
	
	\item[Reila:] Before we go ahead, I'd like to mention another intriguing consistency, namely that $\overleftarrow{9}=-1$ fits perfectly with the geometric series. We have
	$$\ldots 9999 = 9 \cdot \ldots 111 = 9 \cdot \sum_{k=0}^\infty 10^k = 9 \cdot \frac{1}{1-10} = -1,$$
	you see?
	
	\item[Constance:] And again, it's simply wrong, because you ignore the divergence of the series.
	
	\item[Reila:] But shouldn't you be happy because a rule from an old theory generalizes to our new theory?
	
	\item[Constance:] I still wouldn't call it a new theory, but must admit that your observation is at least not a bad omen.

	\item[Reila:] Good girl. Anyway, we found out that $\overleftarrow{9}$ seems to be another representation of $-1$. I wonder if there are other non-natural numbers that can be represented by our new `numbers'.
	
	\item[Nana:] Of course! While you folks were discussing esoteric stuff I found out that all negative integers have such a representation. For instance, $-17$ equals $\overleftarrow{9}83$, because $\overleftarrow{9}83+17=0$. So, any negative integer, say $-x$, is represented by $\overleftarrow{9}$ followed by something. And that something is $100-x$ if $x$ is a two-digit number, $1000-x$, if $x$ a three-digit number, and so on.
	
	\item[Reila:] Well done, Nana. What do we want to consider next? Fractions?
	
	\item[Forest:] That shouldn't be too hard to find out. Let's take one third for a start. We simply have to look for a number with infinitely many digits to the left which gives 1 when multiplied by 3. The last digit, say $a_1$, has to satisfy $3a_1 \equiv 1 \pmod{10}$. The only possible choice is $a_1=7$, and as $3\cdot 7 =21$ this gives a carry of $2$. Therefore, $3a_2+2 \equiv 0 \pmod{10}$, and we obtain $a_2=6$. Again, we have the carry $2$, and obtain $a_3=6$, and so on. Thus, $1/3 = \overleftarrow{6}7$.
	
	\item[Reila:] Not bad, thereby we also have $-1/3=\overleftarrow{3}$, which is certainly intriguing. Let us determine $1/7$. We have $7a_1 \equiv 1 \pmod{10}$, so $a_1=3$. Next we have $7a_2 +2 \equiv 0 \pmod{10}$, so $a_2=4$ and the next carry is $3$. Thus, $7a_3+3 \equiv 0 \pmod{10}$, which implies $a_3=1$.
	
\end{description}

\textit{Reila is scribbling and murmuring for some time...}

\begin{description}
	
	\item[Reila:] That's what I thought: We have $1/7=\overleftarrow{285714}3$ and therefore, $-1/7=\overleftarrow{142857}$.
	
	\item[Preston:] And that's just the same block of repeating digits as in the ordinary decimal representation of $1/7$, namely $0.\overline{142857}$.
	
	\item[Reila:] There is a simple reason for that phenomenon. If you multiply $0.\overline{285714}$ by $7$, you get $0,999....$ and this equals $1$, which explains the ordinary decimal expansion of $1/7$. But that means also that you get $...999$ if you multiply $\overline{285714}$ by $7$. And as we believe in $\overleftarrow{9}=-1$, we obtain that $-1/7$ has the same block of digits repeatedly to the left as the decimal representation of $1/7$ has to the right.
	
	\item[Forest:] Well done, Reila. And this easily generalizes, as it not hinges on the numerator to be 1. I claim that $\overleftarrow{a_1\ldots a_k} = -0,\overline{a_1\ldots a_k}$.
	
	\item[Reila:] Can you prove it?
	
	\item[Forest:] Sure. We have $0,\overline{a_1\ldots a_k} = \frac{a_1\ldots a_k}{10^k-1}$, where $a_1\ldots a_k$ is meant to be the number with the digits $a_1,\ldots a_k$ in decimal representation and not a product or so.
	
	\item[Reila:] Of course, go ahead.
	
	\item[Forest:] Therefore, $-0,\overline{a_1\ldots a_k} = \frac{a_1\ldots a_k}{1-10^k}$. Thus, to prove my claim, I simply have to show that $\overleftarrow{a_1\ldots a_k} \cdot (1-10^k) = a_1\ldots a_k$. But this is obvious since $\overleftarrow{a_1\ldots a_k} \cdot 10^k = \overleftarrow{a_1\ldots a_k} \underbrace{0\ldots 0}_{k}$.
	
	\item[Reila:] Well deduced.

	\item[Forest:] Thanks. I'd like to summarize what we have achieved so far. Formally we consider the set $L=\{(a_n)_{n \in \mathbb{N}} ~|~ a_n \in \{0,1,\ldots,9\} \textnormal{ for all } n\}$.\footnote{The letter $L$ is taken from another interesting article on $10$-adic number representation, namely \cite{Rich} and comes from `Leftist numbers'. Note that in \cite{Rich} elements of $L$ can have finitely many digits after the decimal point (in the further course we will also allow this).} For simplicity, we denote $a \in L$ without parentheses and from right to left as $a=\ldots a_4a_3a_2a_1$. We define addition and multiplication as in the naturals, where carries can go on `infinitely' to the left. The naturals can be identified as a subset of $L$, by noting that for a finite digit string $a_k\ldots a_1$ we can write $\overleftarrow{0}a_k\ldots a_1$; in particular we have $0=\overleftarrow{0} \in L$. Each negative integer can be represented by an element of $L$. Thus, $(L,+,\cdot)$ is a commutative unitary ring.
	
	\item[Preston:] It's so easy to make you happy. But an algebraic structure isn't worth anything in its own right if it has no use.
	
	\item[Forest:] And that's simply not true. Haven't we just mentioned non-euclidean geometries? How many times in history has science benefited greatly from  theories that have been developed and studied previously by mathematicians in their own right?  
	
	\item[Reila:] Calm down, boys. We haven't wasted too much effort on our new theory until now, have we? What about inverse elements of multiplication?
	
\end{description}

\section{Exploring the new area}\label{sec:exploring}

\begin{description}
	
	\item[Forest:] We have already seen that there is an inverse element of $7$, but I think we get problems in general. Take $2$ for instance. There can't be an element $a=\ldots a_2a_2a_1 \in L$ with $2 a = 1$ because this would imply $2a_1 \equiv 1 \pmod{10}$, which cannot be.
	
	\item[Reila:] We could allow finitely many digits right of the decimal point, then $1/2=\ldots 000.5$ would have its ordinary representation.
	
	\item[Constance:] Could we please stop constantly changing the prerequisites we agree on!
	
	\item[Forest:] Constance is right, let's stick to the set $L$ for the time being.
	
	\item[Preston:] I think among the rationals only pure recurring decimals can be represented by elements of $L$.
	
	\item[Nana:] Why's that?
	
	\item[Preston:] It's the same argument as for $1/2$. Considering a rational decimal that is not pure repeating means considering a reduced fraction $m/n$ with $\mbox{gcd}(n,10)>1$. Now, let's assume there was an element $a=\ldots a_3a_2a_1 \in L$ with $m/n=a$. If $n$ was a multiple of $2$ or $5$, we'd have $m= n \cdot m/n = n \cdot \ldots a_3a_2a_1$. Thus, the last digit of $m$ would be a multiple of $2$ or $5$, which contradicts that $m/n$ is reduced.
	
	\item[Forest:] But we could still get a field if we allowed finitely many digits right of the decimal point...
	
	\item[Constance:] \texttt{grunts}
	
	\item[Reila:] We keep that in mind, Forest. Right now, I have another question. We have seen different `new' factorizations of $1$, namely $3 \cdot \ldots 6667 = 1$ or $7 \cdot \overline{142857} = 1$. I wonder if we can also find `new' factorizations of other numbers...
	
	\item[Constance:] Oh, right. How about, for example, $0$?
	
	\item[Nana:] You are just being silly. You can't multiply two non-zero numbers and get zero!
	
	\item[Reila:] Let's see, Nana. Indeed, that is a good question: Can we factor $0$ in some non-trivial way? In other words, are there zero divisors in $L$?
	
	\item[Preston:] Ok, let's start by assuming we have $a=\ldots a_3a_2a_1,b= \ldots b_3b_2b_1 \in L$ with $ab = 0$.
	
	\item[Nana:] You can't just assume that! That is the weird question we are trying to answer!
	
	\item[Preston:] You're right, Nana. But the point of this assumption is not to take an answer to this question for granted, but to get the investigation going: My plan is to derive as much information as I can on these fictitious $a$ and $b$; maybe we get enough information to actually find an example. And once we have that, we can forget the assumption, because the question is answered.
	
	\item[Forest:] Of course, it could also happen that we just run into contradictions, and then we have seen that such $a$ and $b$ do not exist, as you suspect, Nana.
	
	\item[Preston:] Right. Anyway, let's start. We obtain immediately that $a_1 b_1 \equiv 0 \pmod{10}$. We want $a$ and $b$ to be different from $0$, so let us ensure this by stipulating $a_{1}\neq 0$ and $b_{1}\neq 0$. Then we must have $a_1=2$ and $b_1=5$. Or vice versa, but let's continue with this choice of $a_1,b_1$. Following the ordinary multiplication algorithm this gives a carry of $1$, so we can conclude that $a_1b_2 + a_2b_1+1 \equiv 0 \pmod{10}$ in the next step. In other words, $5a_2+2b_1 \equiv 9 \pmod{10}$.
	
	\item[Reila:] This congruence has several solutions, any odd digit as $a_2$ combined with $b_2 \in \{2,7\}$ will do the job.
	
	\item[Preston:] So let's make a choice, say $a_2=1, b_2=2$, and see if we can find $a_3$ and $b_3$.
	
	\item[Constance:] Hmm..., let's see if we can simplify this. Until now, we did two steps following the multiplication algorithm to find zero divisors. More general -- and regardless of the algorithm -- in step $n$ we want to find two $n$-digit numbers $a, b$ with $ab \equiv 0 \pmod{10^n}$. It's clearly sufficient that $a \equiv 0 \pmod{2^n}$ and $b\equiv 0 \pmod{5^n}$.
	
	\item[Reila:] I see what you are getting at, Constance. We could do this inductively. Assume that we have found $a = a_n\ldots a_1$ and $b=b_n \ldots b_1$ with $a \equiv 0 \pmod{2^n}$ and $b \equiv 0 \pmod{5^n}$. How can we find $a_{n+1}$ and $b_{n+1}$ such that $10^n a_{n+1} + a \equiv 0 \pmod{2^{n+1}}$ and $10^n b_{n+1} + b \equiv 0 \pmod{5^{n+1}}$?
	
	\item[Preston:] For $a_{n+1}$ it's pretty easy. Write $a=2^n a'$, so $10^n a_{n+1} + a = 2^n(5^n a_{n+1} + a')$. As we only have to ensure that the expression in parentheses is even, we can choose any even digit as $a_{n+1}$ if $a'$ is even, and any odd digit otherwise.
	
	\item[Reila:] Well done, Preston. That coincides perfectly with our earlier result and generalizes easily to $b_{n+1}$. Write $b=5^n b'$, so $10^n b_{n+1} + b = 5^n(2^n b_{n+1} + b')$ and we have to choose $b_{n+1}$ such that the expression in parentheses is congruent to $0$ modulo $5$. As $3$ is inverse to $2$ modulo $5$ this means that we have to choose $b_{n+1}$ according to $b_{n+1} \equiv 3^n(-b') \pmod{5}$. And of course there are always two such digits.
	
	\item[Constance:] Zero divisors galore! In particular this means that there is not even an ordering of your domain of calculation that is compatible with multiplication.
	
	\item[Nana:] Why not?
	
	\item[Constance:] Consider a product $ab=0$ with $a,b \ne 0$. Then $a,b$ have to be positive or negative. If both were positive the product would have to be positive as well, but it is zero. Similar problems arise whatever `signs' $a$ and $b$ have.
	
	\item[Reila:] We started looking for the largest natural number and ended up discovering a domain where `larger' doesn't even make sense.\footnote{This was actually pointed out by one of the participants of an enrichment course in which the authors treated 10-adic numbers with school students.} Isn't it fascinating?
	
	\item[Constance:] Fascinating? It only shows how pointless your funny `new numbers' are.  
	
	\item[Forest:] Constance, I find it really strange how you could think that something is pointless because it is new. The existence of zero divisors is just another    property of our new numbers, not a reason to reject them. And concerning the ordering, I reply with a quote from Lord of the Rings: `We have not found what we sought, but what have we found?'\footnote{\cite{Tolkien}, p. 191}

	\item[Constance:] I understand that this is best defended by a quote from a fantasy story, just as in the case of your non-associative `addition'. Fortunately, I have the voice of innocence on my side in both cases: Nana made it quite clear that she intended her `numbers' to have neither violations of non-associativity of addition nor zero divisors. I mean, just look at this mess! We started out searching the `largest natural number', we found a candidate, which, however, turned out to be not quite so large, but negative. We changed the rules numerous times to make even something as basic as associativity of addition work, and as soon as we introduce multiplication, the next disaster happens. What does it require for you to accept that there is just nothing there?
	
	\item[Forest:] A contradiction. Which is not just something that contradicts your or Nana's expectations, but an internal contradiction. Things can be arbitrarily strange, but they cannot contradict themselves. 
	
	\item[Constance:] Seriously, Forest? That is all that is required to make you happy? I can tell all kinds of crazy stories and you will just listen as long as I don't contradict myself?
	
	\item[Forest:] You should not underestimate the criterion of consistency. Staying consistent is quite hard, just think of Frege. 
	
	\item[Constance:] A crazy story remains a crazy story, even if it is hard to tell it. But I wonder: Why do you even insist on consistency? Given your focus on `manipulation rules', I would expect you to go the final mile and simply accept inconsistency as `just another property'?\footnote{For this question and the following discussion see e.g. \cite{Becker}, p.~33 (free translation): \textit{`Let us repeat the question: Of what significance is the stipulation of consistency in Hilbert's theory? [...] It can no longer be its serving as a conditio sine qua non for the truth! The answer is fairly curious; it reads: Although the consistency of Hilbert's formal mathematics is not an indispensable requirement for the truth of the theorems, it is still the conditio sine qua non for the continuation of the process of deduction. Thus, by consistency, and by consistency alone, the mathematical `formula game' is protected from a premature abort.'}
	}
	
	\item[Forest:] Don't be silly, Constance. From an inconsistent system, we can easily derive anything. They are just trivial and therefore provide no worthwhile subject of research.
	
	\item[Constance:] Aha! So your ultimate criterion turns out to be hedonistic: You accept that an object exists if it is in your interest to do so. 
	
	\item[Forest:] Not quite. Math is the study of formal systems. You won't hear me talking about `existence' of mathematical objects in any sense that transcends derivability in the formal system we  happen to be currently working in. Inconsistent systems are, indeed, just formal systems with a certain property. My `hedonism', as you call it, only comes in when I decide not to work with them. Certainly, you will not claim that my decision to work on one mathematical topic rather than another would be part of mathematics. You, on the other hand, seem to have a much stronger concept of existence, and the single `criterion', if you can call it that, I have seen so far from you is that you only accept that $X$ exists if you learned about $X$ more than $10$ years ago.
	
	\item[Constance:] It is not the worst criterion that something has withstood the test of time, and certainly more than just subjective interest. Ultimately, I don't really care whether the reasons for rejecting nonsense are considered to be part of mathematics or not. Well-established theories prevail and flourish, while fashionable nonsense is being proposed, gathers a brief interest for the sake of curiosity and is then put aside. Just look at your non-associative addition: No one else but you was willing to go ahead with that. The topic was dropped, just as these new `numbers' with their zero divisors will be dropped and forgotten.  Right, Reila?
	
	\item[Reila:] Not quite. It is indeed a bit much to accept a non-associative operation as a generalization of addition, but structures with zero divisors are quite sensible. Just think of residue rings modulo non-primes or matrix rings. The question is, does that lead us anywhere?
	
	\item[Preston:] The real question is, is it good for anything? I am quite willing to accept very strange and even inconsistent theories when they can be put to good use. If it can, it would be crazy not to use it. It will be used, we will get used to it, and with every year, it will sound less crazy. The expositions will become more and more down-to-earth, more and more people will know about it, and after a few decades, folks like Constance will scoff at any criticism of such a well-established theory. On the other hand, if it cannot be applied, it is a pointless discussion who has the better fantasy universe.
	
	\item[Forest:] We have to leave it for another time whether or not `applicability' is a good criterion for the quality of scientific theories. As a criterion for or against investigating a certain topic, it is no great help: We can only apply what we can handle, and learning how to handle something is just what the investigation is about. I propose to continue research rather than wasting our time on meta-discussions. Does anyone have a question for us?

	\item[Reila:] I do. We have learned that some rationals are not in $L$ and some are. I wonder if  $L$ contains any irrational numbers. What do you think?
	
	\item[Constance:] I think it contains nothing but thoroughly irrational numbers, except for the few naturals you kindly offered a place in your realm of madness.
	
	\item[Reila:] Any constructive comments?
	
	\item[Nana:] Well, you can't just have whatever you like. For instance, $\sqrt{2}$ cannot be in $L$. If $a\in L$ ends with an even digit, its square will be divisible by 4, while if it ends with an odd digit, the square will be odd. So $a^2=2$ is not possible.
	
	\item[Preston:] That's right. What about $\sqrt{5}$? Let's assume we have $a=\ldots a_3a_2a_1 \in L$ with $a^2=5$. From $a_1^2 \equiv 5 \pmod{10}$ we can conclude $a_1=5$.
	
	\item[Forest:] So we need $a_2$. If we want $(\ldots a_3a_2a_1)^2=\ldots 005$ we should have $(10a_2+a_1)^2 \equiv 5 \pmod{100}$. Hmm... that doesn't look too good. We have $(10a_2+a_1)^2 = 100a_2^2+ 20a_1a_2 + a_1^2 \equiv 25 \pmod{100}$. That was bad luck, I guess. Let's check another square root, say $\sqrt{11}$...
	
	\item[Constance:] Hang on! I still distrust this whole business, but just trying out random numbers will certainly not get us anywhere. We should conceptualize this a little bit. Let $q \in \mathbb{N}$ be no perfect square. For $\sqrt{q}$ to be in $L$ it is apparently necessary that $q$ is a quadratic residue modulo $10^n$ for every $n$, because in each step we have to find a $n$-digit natural number $x$ with $x^2 \equiv q \pmod{10^n}$.

	\item[Forest:] Okay, let's see. By the Chinese Remainder Theorem\footnote{see e.g. \cite{Rosen}, pp.~293-295}, $q$ is a quadratic residue modulo $10^n$ if and only if it is a quadratic residue modulo $2^n$ and modulo $5^n$.
	
	\item[Reila:] But isn't that just a reformulation of the problem? We still have to find solutions to these congruences for all natural numbers.
	
	\item[Preston:] You're wrong, Reila. There's a very handy result from number theory that does the job now, namely Hensel's Lemma\footnote{For a good presentation of Hensel's Lemma and the mentioned corollary see \cite{Hill}, pp.101-116.}. It immediately implies that if there is a positive integer $k$ such that $f(x)\equiv 0 \pmod{p^k}$ has a solution $x$ that satisfies $f'(x) \not\equiv 0 \pmod{p}$, then $f(x) \equiv 0 \pmod{p^n}$ has a solution for every positive integer $n$. Here, $f$ should be a polynomial with integer coefficients and $p$ prime.
	
	\item[Reila:] I see. We apply this to $f(x)=x^2-q$. So if for instance $x^2 \equiv q \pmod{5}$ has a solution $x$ with $2x \not\equiv 0 \pmod{5}$ we can conclude that $q$ is a quadratic residue modulo $5^n$ for every $n$.
	
	\item[Constance:] Not bad. Let us for the time being stipulate that $q \not\equiv 0 \pmod{5}$. Then we obtain that $q$ is a quadratic residue modulo $n$ for all $n$ if and only if  $q \equiv \pm 1 \pmod {5}$.
	
	\item[Reila:] Fine, but what about the congruences modulo powers of $2$? As we always have $2x \equiv 0 \pmod 2$, we can't use Hensel's Lemma for them. 
	
	\item[Preston:] I remember another Corollary\footnote{For this corollary see also \cite{Hill}, pp.101-116.} that follows from Hensel's Lemma, a variant that is tailor-made for square roots of integers modulo powers of $2$. It says for odd $q$ that if $x^2 \equiv q \pmod{8}$ has a solution, then $x^2 \equiv q \pmod{2^n}$ has a solution for every $n$.
	
	\item[Reila:] That's nice! Thus, if we also stipulate $q \not\equiv 0 \pmod{2}$, by remembering that odd squares are always congruent to $1$ modulo $8$ we obtain that $q$ is a quadratic residue modulo $2^n$ for all $n$ if and only if $q \equiv 1 \pmod{8}$.
	
	\item[Constance:] OK, let me summarize: We've found out that a natural number $q$ that is no perfect square and that satisfies $\gcd(q,10)=1$ has a square root in $L$ if and only if $q \equiv 1 \pmod{40}$ or $q \equiv 9 \pmod{40}$.

	\item[Forest:] I'm not sure about the if part.
	
	\item[Reila:] Why not?
	
	\item[Forest:] Having a finite sequence of length $k$ for every $k$ does not mean to have one infinite sequence.
	
	\item[Nana:] I don't understand a word!
	
	\item[Forest:] Ok, let me explain. We have a quite interesting issue here, and I'm sure you can understand the problem. We are interested in sequences with a certain property, and Reila says that we can find such sequences of arbitrary but finite length. But what we need is to find an infinite sequence with that property.
	
	\item[Nana:] Isn't that the same? 
	
	\item[Constance:] For every $k$ there is a strictly increasing sequence of length $k$ of natural numbers that is bounded. Do you conclude from this fact that there is an infinite strictly increasing sequence of natural numbers that is bounded?
	
	\item[Nana:] Ah, I see the difference now!
	
	\item[Reila:] (stares into space, deep in thought)\\
	
	\textit{Prof. König comes around the corner and beckons.}\footnote{K\"onig's Lemma states that an infinite, finitely branching tree must have an infinite branch. It is an ubiquitous tool in many areas of mathematics, though it is rarely made explicit; see \cite{Simpson}. In this case, our tree would consist of finite sequences of digits $s$ such that $s^{2}\equiv q$ (mod $10^{|s|}$), where $|s|$ denotes the length of $s$ and we freely confuse digit sequences with the decimal numbers they represent.}\\
	
	\item[Reila:] Now I got it! For every $k$ there is a sequence $s_k$ of digits of length $k$ with the desired property, namely that $s_k^2 \equiv q \pmod{10^k}$. Each of these has a first digit and there are only finitely many first digits, so there is one digit $d_1$ that occurs infinitely often. There are infinitely many sequences of length $k$ with the desired property having the first digit $d_1$, and each of these has a second digit. So there is one second digit $d_2$ that occurs infinitely often. Continuing in this way we construct an infinite sequence of digits $d_i$, which will have the desired property.
	
	\item[Constance:] Well done! And what about numbers $q$ with $\gcd(q,10)\ne 1$?
	
	\item[Preston:] I used the time you spent on figuring out the last step to think about that. Essentially we already have the right condition, see?\\
	
	\textit{Preston hands a sheet of paper to the others. After a short while of silent reading everybody nods in agreement.}\footnote{A natural number $q$ that is no perfect square has a square root in $L$ if and only if it is of the form $q=2^{2\ell}5^{2k} q'$, where $\mbox{gcd}(q',10)=1$, $k, \ell \in \mathbb{N}_0$, and $q' \equiv 1 \pmod{40}$ or $q' \equiv 9 \pmod{40}$. This follows from the above seen arguments and the fact that natural numbers that have a square root in $L$ cannot contain odd powers of $2$ or $5$ in their prime factorization. The latter statement is a generalization of the fact that that 2 cannot have a square root in L (see Nana's argument on the bottom of page 13) and can be derived similarly.}\\
	
	\item[Preston:] Well done! I wonder what we can do about roots of higher degree. 
	
	\item[Reila:] Well, for third roots, things just get easier because we can simply apply Hensel's lemma. For simplicity, let's again stipulate that $\text{gcd}(q,10)=1$. Now, for some such natural number $q$, we want to solve $x^{3}-q=0$ modulo $2$ and modulo $5$. Say $f(x)=x^{3}-q$. If these congruences are solvable, we then only need to check that a solution $x$ satisfies $f^{\prime}(x)\neq 0$, i.e., $3x^{2}\neq 0$ modulo $2$ and $5$, but modulo $2$ and $5$, this is just saying that $\text{gcd}(x,10)=1$, which is clear anyway. 
	
	\item[Preston:] For example, we have $1^3\equiv 3$ (mod $2$) and $2^3\equiv 3$ (mod $5$), so we know that there must be a third root of $3$.
	
	\item[Reila:] Actually, any odd number has a third root modulo $2$, and, since every residue modulo $5$ is a third power ...
	
	\item[Constance:] ... since $\text{gcd}(3,\phi(5))=\text{gcd}(3,4)=1$ and by the Chinese remainder theorem...
	
	\item[Reila:] ... we obtain third roots for every natural number that has no common divisor with $10$. And, by Constance's remark, the same will work for any odd exponent that is not divisible by $5$. But what about, for example $4$ and $5$?
	
	\item[Preston:] I am sure they can be brought home with quite similar techniques. But now that I have suspended my interest in applications for quite a while due to Forest's call for patience, I wonder: After quite a bit of investigation of these new numbers, are they actually good for anything?
	
	\item[Constance:] Good point, Preston. I let myself be persuaded to go along with the investigation up to here, but I still think that it was only worth the effort if this new form of number representation has some useful implications in usual fields of mathematics.
	
	\item[Forest:] Although I don't agree with you, I'm not opposing our moving the investigation in that direction now. Preston, have you already something in mind?

\end{description}

\section{Applying the new tools}

\begin{description}
	
	\item[Preston:] The whole thing looks very similar to the method of complements used in computer science.
	\item[Constance:] Our number representations reminded me of that, too. The difference is that the method of complements comes with the restriction to a finite subset of the rationals. With our number representation we don't have this restriction anymore. We can represent all rationals.
	\item[Forest:] That's not true for the set $L$ we agreed on, whose members are supposed to have no digits after the decimal point. We've discussed that we cannot represent any reduced fraction whose denominator contains one of the prime factors 2 or 5.
	\item[Preston:] That's true, so to get our hands on, let's remove this constraint. I suggest to consider all numbers that can be represented with infinitely many digits, but only finitely many may occur after the decimal point.
	\item[Forest:] Alright, but there's another problem. As these representations have infinite lengths, we cannot use them for computer arithmetic, unless we truncate them at some point. So we haven't gained anything, have we?
	\item[Preston:] I see the problem, but there's an easy solution. If we assume in addition that the digits to the left are repeating from some point on, each representation can be written down in finite length by using Reila's arrow notation, and still the set of numbers we can represent contains at least all rationals.
	\item[Reila:] To be more precise, it's exactly the rationals that we consider then.
	\item[Forest:] Sure?
	\item[Reila:] Consider such a representation, say $\overleftarrow{a_1\ldots a_k}b_1\ldots b_l.c.$ This equals
	$$10^k \overleftarrow{a_1\ldots a_k} + b_1\ldots b_l + 0.c_1\ldots c_r.$$
	And as we already know that $\overleftarrow{a_1\ldots a_k} = -0.\overline{a_1\ldots a_k}$ is rational, all three summands are rational.\\
	
	\textit{Everybody nods in agreement.}\\
	
	\item[Preston:] Okay, so we have a new number representation to do computer arithmetic on the rationals.
	\item[Constance:] Fine, but what are the advantages of our new number representation in comparison to the usual decimal representation with regards to computer arithmetic?
	\item[Preston:] Well, there are some disadvantages with the usual decimal representation, even with the basic arithmetic operations. Take for instance simple addition. The usual algorithm proceeds from right to left...
	\item[Reila:] Where's the problem?
	\item[Preston:] Add $2.506218...$ and $1.493781...$!
	\item[Reila:] The result starts with $3.999999$, unless ..., okay, I see the problem. If there was a carry, the result would begin with $4.000000$. So I can't tell you even a single digit of the result.
	\item[Forest:] But that's only because Preston didn't define the summands properly.
	\item[Preston:] Agreed. But it's still a flaw. Although I specified the summands up to some small error, you cannot tell any digit of the result. And the fact the number representation is not unique, like with $3.\overline{9} = 4$ is also an annoyance.
	\item[Constance:] You clearly have a point, but it doesn't knock my socks off.
	\item[Preston:]Another thing is the problem with negative summands.
	\item[Reila:] What's wrong with them?
	\item[Constance:] I can imagine what Preston's getting at. If you want to perform an addition and the summands might be negative, you cannot simply add. You first have to check, which of the summands are negative. If both are non-negative, you know that the result is also non-negative and you add. If both are negative, you know that the result is negative, and you add the absolute values to obtain the absolute value of the result.
	\item[Preston:] Exactly. And the worst case is when one is negative and the other is positive, because then you first have to check which one you have to subtract from the other one by comparing the absolute values.
	\item[Reila:] I see. Take for instance $19+(-17)$ versus $2+(-17)$. Traditionally, in the first case we have to subtract  $17$ from $19$ to get $2$. And in the second case we have to subtract $2$ from $17$ and, because we know that the result has to be negative, we get $-15$. What we actually have to do is quite complicated and depends on the numbers, although both calculations are simple additions.
	\item[Nana:] And with my numbers and Reila's arrow notation we don't have any trouble, ... 
	\item[Forest:] Nana, at the moment we are considering rational numbers. So, the numbers are the same, it's just that we use a different representation for them.
	\item[Reila:] But without Nana's input we wouldn't be considering anything new! Anyway, I think what Nana was getting at is that we have no trouble with signs, because we can always simply add. In the first case the task reads $19+\overleftarrow{9}83$ and the usual addition algorithm can be applied.
	\item[Nana:] It's easy, see?
	\begin{small}
		$$\begin{array}{r}
			\ldots 00019\\
			+\phantom{0}\ldots{99983}\\[-2mm]
			\rule{13mm}{0.4pt}\\
			\ldots 00002
		\end{array}$$
	\end{small}
	The second task reads $2+\overleftarrow{9}83$, and we can deal with it in the exact same way
	\begin{small}
		$$\begin{array}{r}
			\ldots 00002\\
			+\phantom{0}\ldots{99983}\\[-2mm]
			\rule{13mm}{0.4pt}\\
			\ldots 99985
		\end{array}$$
	\end{small}
	\item[Preston:] Not bad, eh?
	\item[Reila:] Agreed. And subtraction is thereby also no issue. What about multiplication?
	\item[Preston:] Should work as usual, with the advantage of not having to think about signs. Normally, we would have to inspect the signs of the factors to decide whether the result will be positive or negative, keep them in mind, do the multiplication with the absolute values, and amend the result with the correct sign afterwards. In our representation we can simply multiply.
	\item[Nana:] I wanna try it, say $\overleftarrow{37}14 \cdot 23$. Let's see...\\
	
	\textit{Nana is scribbling and murmuring for a while.}\\
	
	\item[Nana:] No problem ...
	\begin{small}
		$$\begin{array}{r}
			\ldots 7373714\cdot 23\\[-2mm]
			\rule{22mm}{0.4pt}\\
			\ldots 47474280\\
			\ldots 12121142\\[-2mm]
			\rule{22mm}{0.4pt}\\
			\ldots 59595422
		\end{array}$$
	\end{small}
	\item[Reila:] Alright, but more interesting would be an example where both factors have infinitely many digits.
	
	\item[Nana:] Here goes, it's fun for me. Let's do $\overleftarrow{3}5\cdot \overleftarrow{8}3$. We have $\overleftarrow{3}5\cdot 3 = 5$ and $\overleftarrow{3}5 \cdot 8 = \overleftarrow{6}80$. That's all we need.\\
	
	\textit{Nana is scribbling and murmuring for a while.}\\
	
	\item[Nana:] An here's what I've got.
	\begin{small}
		$$\begin{array}{l}
			\ldots 000000000000\bm{5}\\
			\ldots 666666666680\\[-2mm]
			\rule{26mm}{0.4pt}\\
			\ldots 66666666668\bm{0}\\
			\ldots 66666666680\\[-2mm]
			\rule{26mm}{0.4pt}\\
			\ldots 3333333334\bm{8}\\
			\ldots 6666666680\\[-2mm]
			\rule{26mm}{0.4pt}\\
			\ldots 000000001\bm{4}\\
			\ldots 666666680\\[-2mm]
			\rule{26mm}{0.4pt}\\
			\ldots 66666668\bm{1}\\
			\ldots 66666680\\[-2mm]
			\rule{26mm}{0.4pt}\\
			\ldots 3333334\bm{8}
		\end{array}$$
	\end{small}
	\item[Nana:] And from this point on it's repeating, so we have $\overleftarrow{3}5\cdot \overleftarrow{8}3 = \overleftarrow{148}05$. Do you understand my scheme? The first line is the result of $\overleftarrow{3}5\cdot 3$, and its last digit is the last digit of the result. The line below is the result of $\overleftarrow{3}5 \cdot 8$. It's indented, because actually we did $\overleftarrow{3}5 \cdot 80$ here. I sum up, and the last digit of the sum is the second last digit of the result, and so on.
	\item[Reila:] So, the multiplication algorithm works fine. What about ...
	\item[Preston:] ... division, sure. I think it's most interesting. And if I'm not mistaken, it'll turn out to be quite nice.
	\item[Constance:] Why's that?
	\item[Preston:] Because we have to proceed from right to left. So with our number representation all four algorithms work in the same direction. And there's another advantage, but you will see. Nana, would you like to do $\overleftarrow{5}8:13$?
	\item[Nana:] Sure, let's see. As we have to start from the right, we are first looking for the last digit of the result. But that's easy, we just have to find a digit $d$, such that $d\cdot 13$ ends in 8. Thus, $d=6$ is the last digit of the result.
	\item[Preston:] Exactly. that's how it works. Go ahead.
	\item[Nana:] We subtract $6\cdot 13 = 78$ from $\overleftarrow{5}8$, that gives  $\overleftarrow{5}480$. And now we think about the second last digit of the result. Of course, we have to ignore the $0$ at the end. So, we are again looking for a digit $d$ with $d\cdot 13 \equiv 8 \pmod{10}$. Thus, the second last digit of the result is also $6$.
	\item[Reila:] Now,  $\overleftarrow{5}48-6\cdot 13 =  \overleftarrow{5}48 - 78 =  \overleftarrow{5}470$. And $d\cdot 13 \equiv 7 \pmod{10}$ holds if and only if $d=9$. So the results end with $966$.
	\item[Nana:] Let me put it into a nice scheme again.\\
	
	\textit{Again, Nana is working feverishly.}\\
	
	\item[Nana:] And here you are.
	\begin{small}
		$$\begin{array}{l}
			\phantom{00000000} \ldots 5558:13 = \overleftarrow{18}4966\\
			\phantom{0000000000000} 78\\[-2mm]
			\phantom{0000000000} \rule{8mm}{0.4pt}\\
			\phantom{0000000} \ldots 5548\\
			\phantom{000000000000} 78\\[-2mm]
			\phantom{000000000} \rule{8mm}{0.4pt}\\
			\phantom{000000} \ldots 5547\\
			\phantom{00000000000} 27\\[-2mm]
			\phantom{00000000} \rule{8mm}{0.4pt}\\
			\phantom{00000} \ldots 5552\\
			\phantom{0000000000} 12\\[-2mm]
			\phantom{0000000} \rule{8mm}{0.4pt}\\
			\phantom{0000} \ldots \color{red}{5554}\\
			\phantom{000000000} 24\\[-2mm]
			\phantom{000000} \rule{8mm}{0.4pt}\\
			\phantom{000} \ldots 5553\\
			\phantom{00000000} 13\\[-2mm]
			\phantom{00000} \rule{8mm}{0.4pt}\\
			\phantom{00} \ldots \color{red}{5554}\\
			
		\end{array}$$
	\end{small}
	
	\item[Reila:] Works just fine. But what's this other advantage you mentioned, Preston?
	
	\item[Preston:] Unlike the usual from-left-to-right algorithm, our new method needs no guessing. You know what I mean, the `How often goes this into that' part. With our method we could find every digit of the result by table look-up.
	\item[Constance:] That's indeed nice.
	\item[Forest:] Nice in this particular example. But our algorithm won't work if the divisor has a prime factor 2 or 5.
	\item[Preston:] Well spotted, Forest. This is because then the divisor has a common factor with our base 10. But we can remove the annoying common factor easily. For instance, if we want to divide by 15, we can just multiply both dividend and divisor by 2. Then we have to perform a division by 30, which we accomplish by dividing by 3 and adjusting the decimal point afterwards.
	\item[Forest:] Agreed. On the other hand, we could solve the problem in a different way. And I think I'd prefer that way. Why don't we just change the base? We should have thought about it before. Base 10 only causes problems, remember the zero divisors we encountered. With a prime base, they'd be gone...\\
	
	\textit{Our friends didn't notice that while they were working a small bunch of people had gathered at the open door. The number theorist Andrew, Tom from the physics department, Marc from computer science, and Dayna, a differential-equations specialist, were watching them intently. When Forest turns around and notices them, they apologize for disturbing and hurry off.}\footnote{For a good survey on the applications of $p$-adic numbers in various fields see \cite{Rozikov}, where, among many other things, it is mentioned that $p$-adic numbers are employed in Andrew Wiles' famous proof of Fermat's Last Theorem. For more information on how Wiles' proof uses $p$-adic numbers see \cite{Gouvea1994}. For a more detailed discussion of the usefulness of the here considered 10-adic number representation with respect to computer science see \cite{Hehner}. It is argued there that this \textit{`novel system for representing the rational numbers [...] allows exact arithmetic, and approximate arithmetic under programmer control  [..., and] is superior to existing coding methods because the arithmetic operations take particularly simple, consistent forms. These attributes make the new number representation attractive for use in computer hardware.'} \cite{Hehner}, p.~1}\\
	
	\item[Preston:] That's a very good suggestion, Forest. But considering that we started looking for the largest natural number and ended up accelerating machine arithmetic, don't you think that we have deserved our afternoon?
	
	\item[Forest:] You're right. But tomorrow we should start analyzing our new number representation with a prime base.\footnote{For a good introduction on $p$-adic numbers see \cite{Gouvea1997}.}
	
	\item[Reila:] I'm in. Actually, I still have a few things bugging me: We investigated and even applied this stuff, but we still don't know what it is. For example how should we interpret the fact that intuitively, the sequence $9,99,999, \ldots$ seems to converge to $-1$? I have an idea which I would expect to advance our investigation.\footnote{For the definition of the so called $p$-adic absolute value and applications thereof, see, e.g., \cite{Gouvea1997}.}
	
	\item[Nana:] Maybe I will stop by again. You folks are smart but sometimes a tad narrow-minded, probably could use my help...
	
\end{description}

\end{document}